\documentclass[12pt,a4paper]{article}
\IfFileExists{ajr.sty}{\usepackage{ajr}}{}

\title{Rigorous modelling of nonlocal interactions determines a macroscale advection-diffusion PDE}

\author{Prof A.J. Roberts
\\School of Mathematical Sciences,
University of Adelaide
\\\url{http://orcid.org/0000-0001-8930-1552}}
\date{21 January 2020}

\usepackage{parskip} 
\usepackage{pgfplots}
\pgfplotsset{compat=newest}

\IfFileExists{ajr.sty}{\usepackage{hyperref}
    \hypersetup{pdftex,
            backref=true,
            hyperindex=true,
            colorlinks=true,
            allcolors=teal!60!black,
            bookmarks=true,
            breaklinks=true}
}{}
\usepackage{amsmath,amssymb,eulervm,defns}
\usepackage{listings}
\iffalse\usepackage{mybiblatex}
\else 
    \usepackage{natbib}
    \AtEndDocument{\raggedright\bibliography{bibexport}}
    \RaisedNamesfalse
\fi

\usepackage{mycleveref}

\usepackage{url,doi}

\allowdisplaybreaks
\Vec r\Vec u\Vec U
\Cal B\Cal W\Cal L\Cal A\Cal V
\Bb R \Bb X \Bb Y 
\def\fL{\mathfrak L}
\newcommand{\intYX}{\int_{\YY}\int_{\XX}}
\def\pdf{\textsc{pdf}}
\def\ide{\textsc{ide}}
\newcommand{\rvd}{\rv'}
\newcommand{\pt}[2][X]{\ifx0#2\else\ifx1#2(x-#1)
  \else\frac{(x-#1)^{#2}}{#2!}\fi\fi}
\newcommand{\qt}[2][X]{\ifx0#2\else\ifx1#2(\xi-#1)
  \else\frac{(\xi-#1)^{#2}}{#2!}\fi\fi}
\newcommand{\xey}{\hat x} 

\def\x{\ensuremath\zeta}
\def\u{\ensuremath{\tilde u}}
\def\r{\ensuremath{\tilde r}}

\begin{document}

\maketitle

\begin{abstract}
A slowly-varying or thin-layer multiscale assumption
empowers macroscale understanding of many physical scenarios
from dispersion in pipes and rivers, including beams,
shells, and the modulation of nonlinear waves, to
homogenisation of micro-structures. Here we begin a new exploration of
the scenario where the given physics has non-local
microscale interactions. We rigorously analyse the
dynamics of a basic example of shear dispersion. Near each
cross-section, the dynamics is expressed in the local
moments of the microscale non-local effects. Centre manifold
theory then supports the local modelling of the system's
dynamics with coupling to neighbouring cross-sections as a
non-autonomous forcing. The union over all cross-sections
then provides powerful new support for the existence and
emergence of a macroscale model advection-diffusion \pde\
global in the large, finite-sized, domain. The approach
quantifies the accuracy of macroscale advection-diffusion
approximations, and has the potential to open previously
intractable multiscale issues to new insights.
\end{abstract}

\tableofcontents

\section{Introduction}

This paper introduces a new rigorous approach to the multiscale challenge of systematically modelling by macroscale~\pde{}s the dynamics of microscale, \emph{spatially nonlocal}, systems.
This approach provides a novel quantified error formula.
Previous research using this type of approach rigorously modelled systems that were expressed as \pde{}s on the microscale.
This previous research encompassed both cylindrical multiscale domains \cite[]{Roberts2013a} and more general multiscale domains \cite[]{Roberts2016a, Bunder2018a}.
But recall that \pde{}s are themselves mathematical idealisations of physical processes that typically take place on microscale length scales.
Hence, here we begin to address the challenges arising when the given mathematical model of a system encodes microscale physical interactions over finite microscale lengths.

Physical systems with nonlocal, microscale, spatial interactions arise in many applications.
In neuroscience,  a spatial convolution expresses the excitatory\slash inhibitory effects of a neurone on a nearby neurone, giving rise to nonlocal neural field equations, and ``have been quite successful in explaining various experimental findings'' \cite[e.g.]{Ermentrout2015}.
Models of free crack propagation in brittle materials invoke microscale \emph{nonlocal} stress-strain laws, called peridynamics \cite[e.g.]{Silling2000}: one challenge is to derive the effective mesoscale \pde{}s from the nonlocal laws \cite[e.g.]{Silling2008, Lipton2014}.
Nonlocal dispersal and competition models arise in biology \cite[e.g.]{Omelyan2018, Duncan2017}.
Other examples are non-local cell adhesion models \cite[e.g.]{Buttenschon2020}.
In this introduction we begin by exploring the specific example of a so-called `Zappa' dispersion in a channel (\cref{seczsd}) in which material is transported by finite jumps along the channel, and also is intermittently thoroughly mixed \text{across the channel.}%

\paragraph{General scenario}
Zappa dispersion is a particular case of the following general scenario---a scenario that is the subject of ongoing research.
In generality we consider a field~\(u(x,y,t)\), on a `cylindrical' spatial domain~\(\XX\times\YY\) (where \(\XX\subseteq\RR\) and where~\YY\ denotes the cross-section).
We suppose the field~\(u\) is governed by a given autonomous system in the form
\begin{equation}
\D tu=\intYX k(x,\xi,y,\eta)u(\xi,\eta,t)\,d\xi\,d\eta\,,
\label{eq:ukint}
\end{equation}
where the given kernel~\(k(x,\xi,y,\eta)\) expresses both nonlocal and local physical effects at position~\((x,y)\) from the field at position~\((\xi,\eta)\), both within the cylindrical domain \(\XX\times\YY\).
We allow the kernel to be a generalised function so that local derivatives may be represented by derivatives of the Dirac delta function~\(\delta\): for example, a component~\(\delta'(x-\xi)\delta(y-\eta)\) in the kernel~\(k\) encodes the differential term~\(-\D xu\) in the right-hand side of~\eqref{eq:ukint}.
In general the physical effects encoded in the kernel~\(k\) may be heterogeneous in space.
But, as is common and apart from boundaries, Zappa dispersion is homogeneous in space (translationally invariant) in which case some significant simplifications ensue. 

The nonlocal system~\eqref{eq:ukint} is linear for simplicity, but we invoke the framework of centre manifold theory so the approach should, with future development, apply to nonlinear generalisations as in previous work on such modelling where the system is  expressed as \pde{}s on the microscale \cite[]{Roberts2013a}.

Our aim is to rigorously establish that the emergent dynamics of the nonlocal system~\eqref{eq:ukint} are captured over the 1D spatial domain~\XX\ by a mean\slash averaged\slash coarse\slash macroscale variable~\(U(x,t)\) that satisfies a macroscale, second-order, advection-diffusion \pde\ of the form
\begin{equation}
\D tU\approx A_1\D xU+A_2\DD xU\,,\quad x\in\XX\,,
\label{eq:dcdtA}
\end{equation}
for some derived coefficients~\(A_1\) and~\(A_2\).%
\footnote{Ongoing research aims to generalise the approach here to certify the accuracy of \pde{}s truncated to \(N\)th-order for every~\(N\).}
This macroscale \pde~\eqref{eq:dcdtA} is to model the dynamics of the microscale nonlocal~\eqref{eq:ukint} after transients have decayed exponentially quickly in time, and to the novel quantified error~\eqref{eq:pderemain}.

\section{Zappa shear dispersion}
\label{seczsd}

This section introduces a basic example system (non-dimensional) of nonlocal microscale jumps by a particle (inspired by W.~R. Young, private communication).
\cref{secmkglm} systematically derives an advection-diffusion \pde~\eqref{eq:dcdtA} for the particle that is valid over macroscale space-time.
Consider a particle in a channel \(-1<y<1\), \(\YY=(-1,1)\), and of notionally infinite extent in~\(x\), \(\XX=\RR\).
Let \(u(x,y,t)\) be the probability density function (\pdf) for  the particle's location: equivalently, view~\(u(x,y,t)\) as the concentration of some continuum material.

The `Zappa' dynamics of the particle's \pdf\ is encoded by 
\begin{equation}
\D tu=\Big[\tfrac1{v(y)}\underbrace{\int_{-\infty}^x e^{-(x-\xi)/v(y)}u(\xi,y,t)\,d\xi}_{{=\ e^{-x/v(y)}\star u},\text{ the convolution  \eqref{eq:conv}}}-u\Big]
+\left[\tfrac12\int_{-1}^1 u\,dy-u\right]
\label{eq:muzap}
\end{equation}
for some jump profile \(v(y)>0\)\,---\(v(y)\) is an effective velocity along the channel.
That is, the kernel of the Zappa system is the generalised function
\begin{align}
k(x,\xi,y,\eta)&=\left[\tfrac1{v(y)}e^{-(x-\xi)/v(y)}H(x-\xi) -\delta(x-\xi)\right]\delta(y-\eta)
\nonumber\\&\quad{}
+\left[\tfrac12-\delta(y-\eta)\right]\delta(x-\xi),
\label{eq:zapker}
\end{align}
where \(H(x)\) is the unit step function.
The nonlocal equation~\eqref{eq:muzap} governs the \pdf\ of the particle in Zappa dispersion through the following two physical mechanisms.
\begin{itemize}
\item We suppose that, at exponentially distributed time intervals with mean one, the particle gets `zapped' across the channel (by a burst of intermittent turbulence for example) and lands at any cross channel position~\(y\) with uniform distribution.  
Consequently the Fokker--Planck \pde~\eqref{eq:muzap} for the \pdf\ contains the terms
\(u_t=\big[\tfrac12\int_{-1}^1 u\,dy-u\big]+\cdots\).

\item Further, suppose that, at exponentially distributed time intervals with mean one, the particle jumps in~\(x\) a distance to the right, a distance which is exponentially distributed with some given mean~\(v(y)\).
Consequently the Fokker--Planck \pde~\eqref{eq:muzap} for the \pdf\ contains the terms
\(u_t=\big[\tfrac1{v(y)}e^{-x/v(y)}\star u-u\big]+\cdots\), in terms of the upstream convolution 
\begin{equation}
e^{-x/v(y)}\star u=\int_{-\infty}^x e^{-(x-\xi)/v(y)}u(\xi,y,t)\,d\xi\,.
\label{eq:conv}
\end{equation}
\end{itemize}
We derive the macroscale model that the cross-sectional mean field~\(U(x,t)\) evolves according to an advection-diffusion \pde: \(U_t\approx A_1U_x+A_2U_{xx}\)\,.
The field~\(U(x,t)\) may be viewed as the marginal probability density of the particle being at~\(x\), averaged over the cross-section~\(y\).
Innovatively, we put the macroscale modelling on a rigorous basis that additionally quantifies \text{the error.}

In particular,  say we choose \(v(y):=1-y^2\) then computer algebra (\cref{app}) readily derives that over large space-time scales, and after transients decay roughly like~\(e^{-t}\), from every initial condition the Zappa system~\eqref{eq:muzap} has the quasistationary distribution \cite[e.g.]{Pollett90} 
\begin{subequations}\label{eqs:zappasm}%
\begin{align}
&u(x,y,t)\approx U+(y^2-\tfrac13)\D xU +(2y^4-\tfrac83y^2+\tfrac{22}{45})\DD xU\,,
\label{eqqspdf}
\\&\text{such that}\qquad
\D tU= -\frac23\D xU+\frac{28}{45}\DD xU+\rho\,,
\quad\label{eq:zappact}
\end{align}
in terms of a macroscale variable here chosen to be the cross-sectional mean,
\begin{equation}
U(x,t):=\frac12\int_{-1}^1u(x,y,t)\,dy\,.
\label{eq:cumean}
\end{equation}
The macroscale \pde~\eqref{eq:zappact} is a precise equality because we include the error terms in our analysis to find a precise, albeit complicated, expression for the final error~\(\rho\). 
The remainder error~\(\rho\) in~\eqref{eq:zappact} has the form
\begin{align}
\rho&:=r_0
+\langle Z_0,\cW_{0:} \cB e^{\cB t}\star\rvd\rangle
+\langle Z_0,\cW_{0:} \rvd\rangle
\nonumber\\&\quad{}
-A_1\langle Z_0,\cW_{1:} e^{\cB t}\star\rvd\rangle
-A_2\langle Z_0,\cW_{2:} e^{\cB t}\star\rvd\rangle
\label{eq:pderemain}
\end{align}
\end{subequations}
where here the convolutions are over time, \(f(t)\star g(t)=\int_0^tf(t-s)g(s)\,ds\)\,, and other symbols are introduced in the next \cref{secmkglm}.
We anticipate this error~\(\rho\) is \begin{itemize}
\item `small' in regions of slow variations in space, small gradients, and 
\item `large' in regions of relatively large gradients such as spatial boundary layers.
\end{itemize}

Then, simply, the macroscale \pde\ model~\cref{eq:zappact} is valid whenever and wherever the error~\(\rho\) is small enough for the application purposes at hand.
The next section includes deriving this error term and clarifies the notation.

\section{Many kernels generate local models}
\label{secmkglm}

Inspired by earlier research \cite[Proposition~1]{Roberts2013a}, this section's aim is to rigorously derive and justify the model~\eqref{eqs:zappasm} that governs the emergent macroscale evolution of Zappa dispersion.
The algebra starts to `explode'---\cref{seccacwm} discusses how to compactly do the algebra in physically meaningful forms, and connect to other mathematical methodologies.

To derive the advection-diffusion model~\eqref{eq:zappact} we truncate the analysis to second order quadratic terms.
Higher-orders appear to be similar in nature, but much more involved algebraically, and are left for later development.

\subsection{Rewrite the equations for local dynamics}
\label{secreld}

Let's analyse the dynamics in the spatial locale about a generic longitudinal cross-section \(X\in\XX\). 
Then invoke Lagrange's Remainder Theorem---which empowers us to track errors---to expand the \pdf\ as
\begin{equation}
u(x,y,t)=u_0(X,y,t)\pt0+u_1(X,y,t)\pt{1}+u_2(X,x,y,t)\pt 2\,,
\label{eq:utrtn}
\end{equation}
where \(u_0:=u\) and \(u_1:=\D xu\) both evaluated at the cross-section \(x=X\), and 
where \(u_2:=\DD xu\) evaluated at some point \(x=\xey(X,x,y,t)\) which is some definite (but usually unknown) function of cross-section~\(X\), longitudinal position~\(x\), cross-section position~\(y\), and time~\(t\). 
By the Lagrange Remainder Theorem, the location~\(\xey\) satisfies \(X \lessgtr\xey \lessgtr x\).
The function~\(\xey\) is implicit in our analysis because it is hidden in the dependency upon~\(x\) of the second derivative~\(u_2(X,x,y,t)\).

Substitute~\eqref{eq:utrtn} into the Zappa nonlocal equation~\eqref{eq:muzap} to obtain
\RaisedNamesfalse
\begin{align}
&\D t{u_0}\pt0+\D t{u_1}\pt{1}+\D t{u_2}\pt 2
\nonumber\\&
=\int_\YY \left[\int_\XX k(x,\xi,y,\eta)\qt{0}\,d\xi\right]
u_0(X,\eta,t)\,d\eta
\nonumber\\&\quad{}
+\int_\YY \left[\int_\XX k(x,\xi,y,\eta)\qt{1}\,d\xi\right]
u_1(X,\eta,t)\,d\eta
\nonumber\\&\quad{}
+\intYX k(x,\xi,y,\eta)\qt{2}u_2(X,\hat\xi,\eta,t)\,d\xi\,d\eta \,.
\qquad\label{eq:sumie}
\end{align}
The effect at cross-section~\(x\) of the \(n\)th~moment of the kernel at cross-section~\(X\) is summarised in the integrals \(\int_\XX k(x,\xi,y,\eta)\qt{n}\,d\xi\)\,.
So define the local \(n\)th~moment of the kernel to be, for every \(n\geq0\),
\begin{align}
k_n(X,y,\eta)&:=\int_\XX k(X,\xi,y,\eta)\qt{n}\,d\xi
\nonumber\\&
=\big[(-v)^n -\delta_{n0}\big]\delta(y-\eta)
+\big[\tfrac12-\delta(y-\eta)\big]\delta_{0n} 
\label{eq:kern}
\end{align}
upon substituting the Zappa kernel~\eqref{eq:zapker}.
This Zappa problem is homogeneous in~\(x\), as are many problems, and so the kernel moments~\(k_n\) are independent of the cross-section~\(X\) (except near the boundary inlet and outlet).

The last integral term in the local expansion~\eqref{eq:sumie} requires special consideration: apply Lagrange's Remainder Theorem to write \(u_2(X,\xi,\eta,t)=u_2(X,X,\eta,t)+(\xi-X)u_{2x}(X,\hat\xi,\eta,t)\) for some uncertain function~\(\hat\xi(X,\xi,\eta,t)\) that satisfies \(X \lessgtr \hat\xi \lessgtr \xi\) for every~\(\eta,t\), and where \(u_{2x}:=\D x{}\left[u_2(X,x,\eta,t)\right]\).
Then the last term distributes into two:
\begin{align*}
&\intYX k(x,\xi,y,\eta)\qt{2}u_2(X,\hat\xi,\eta,t)\,d\xi\,d\eta
\\&=
\int_\YY \underbrace{\int_\XX k(x,\xi,y,\eta)\qt{2}\,d\xi}_{k_2(X,y,\eta)}
{u_2(X,X,\eta,t)}\,d\eta
\\&\quad{}+
\underbrace{\intYX k(x,\xi,y,\eta)3\qt3 u_{2x}(X,\hat\xi,\eta,t)\,d\xi\,d\eta}_{\text{a remainder, with a third $x$ derivative in }u_{2x}} \,.
\end{align*}
Define \(u_2(X,y,\eta):=u_2(X,X,y,\eta)\) for notational consistency with lower moments---see the definition~\eqref{eq:kern}.

The local equation~\eqref{eq:sumie} is exact everywhere, but is most useful in the vicinity of the cross-section~\(X\), that is, for small~\((x-X)\).
Notionally we want to `equate coefficients' of powers of~\((x-X)\) in~\eqref{eq:sumie}, but to be precise we must {carefully} evaluate \(\lim_{x\to X}\) of various \(x\)-derivatives of~\eqref{eq:sumie}.
For example, let \(x\to X\) in~\eqref{eq:sumie}, then
\begin{eqnarray*}
\D t{u_0}&=&\int_\YY k_0(X,y,\eta) u_0(X,\eta,t)\,d\eta
+\int_\YY k_1(X,y,\eta) u_1(X,\eta,t)\,d\eta
\nonumber\\&&{}
+\int_\YY k_2(X,y,\eta) u_2(X,\eta,t)\,d\eta
\\&&{}
+3\intYX k(X,\xi,y,\eta)\qt3 u_{2x}(X,\hat\xi,\eta,t)\,d\xi\,d\eta
\,.
\end{eqnarray*}
Let's rewrite this conveniently and compactly as the integro-differential equation (\ide)
\begin{equation}
\D t{u_0}=\fL_0 u_0+\fL_1u_1+\fL_2u_2
+r_0\,,
\label{eq:du0dt}
\end{equation}
for $y$-operators defined to be, from the moments~\eqref{eq:kern}, 
\begin{equation}
\fL_n u:=\int_\YY k_n(X,y,\eta)u|_{y=\eta}\,d\eta
=\begin{cases}
\tfrac12\int_{-1}^1u\,dy -u\,,&n=0\,,\\
[-v(y)]^n u\,,&n=1,2,\ldots\,.
\end{cases}
\label{eq:deffLn}
\end{equation}
The \ide~\eqref{eq:du0dt} also has the remainder~\(r_0\) which couples the local dynamics to neighbouring locales via~\(u_{2x}\) and is the \(n=0\) case of
\begin{equation}
r_n(X,y,t):=
3\intYX \Dn xnk\Big|_{x=X}\qt3 u_{2x}(X,\hat\xi,\eta,t)\,d\xi\,d\eta \,.
\label{eq:du0dtr}
\end{equation}
Now we can see how this approach to modelling the spatial dynamics works: 
given that the \(y\)-operators~\eqref{eq:deffLn} are evaluated at~\(X\), the spatially local power series with remainder, in \ide{}s like~\eqref{eq:du0dt}, `pushes' the coupling with neighbouring locales to a higher-order derivative term in~\(r_0\), here third-order via the \(u_{2x}\)~factor.
Hence the local dynamics in~\(u_0,u_1,u_2\) are essentially isolated from all other cross-sections whenever and wherever the coupling~\(r_0\) is small enough for the purposes at hand---here when third derivatives are small---that is, when the solutions are, in space, slowly varying enough.

The previous paragraph obtains the \ide\ for~\(u_0\) by simply taking the limit of~\eqref{eq:sumie} as \(x\to X\).
We straightforwardly and similarly obtain \ide{}s for \(u_1\) and~\(u_2\) by finding the limits of spatial derivatives of~\eqref{eq:sumie}:
\begin{subequations}\label{eqs:dundt}%
\begin{align}
\lim_{x\to X} \D x{\eqref{eq:sumie}} &\implies
\D t{u_1}=\fL_0 u_1+\fL_1u_2+r_1\,;
\\
\lim_{x\to X} \DD x{\eqref{eq:sumie}} &\implies
\D t{u_2}=\fL_0 u_2+r_2\,;
\end{align}
\end{subequations}
for local coupling remainders~\(r_1\) and~\(r_2\) defined by~\eqref{eq:du0dtr}.

\subsection{Local-to-global system modelling theory}
\label{seclgsmt}

This section considers the collection of `local' systems as one `global' (in space~\(X\)) system.
Then theory establishes that the advection-diffusion \pde~\eqref{eq:zappact} arises as a globally valid, macroscale, model \pde.

Denote the vector of coefficients \(\uv(X,y,t):=(u_0,u_1,u_2)\),  and similarly for the local coupling remainder \(\rv(X,y,t):=(r_0,r_1,r_2)\).
Then write the \ide{}s~\eqref{eq:du0dt} and~\eqref{eqs:dundt}, in the form of the `forced' linear system 
\begin{equation}
 \frac{d\uv}{dt}=\underbrace{\begin{bmatrix} \fL_0&\fL_1&\fL_2
\\0&\fL_0&\fL_1\\0&0&\fL_0 \end{bmatrix}}_{\cL}\uv+\rv(X,t).
\label{eq:losys}
\end{equation}
for upper triangular matrix\slash operator~\cL.
The system~\eqref{eq:losys} might appear closed, but it is coupled via the  derivative~\(u_{2x}\), through the `forcing' remainders~\rv, to the dynamics of cross-sections that neighbour~\(X\). 

At each locale~\(X\in\XX\), treat the remainder coupling~\rv\ (third-order) as a perturbation (and if this was a nonlinear problem, then the nonlinearity would also be part of the perturbation).
Thus to a useful approximation the global system satisfies the local linear \ode{}s \(\de t\uv\approx\cL\uv\) for each~\(X\in\XX\).
Hence, the linear operator~\(\cL\) is crucial to modelling the dynamics: all solutions are characterised by the eigenvalues of~\cL.
Since~\cL\ is block triangular, a structure exploited previously \cite[\S2]{Roberts2013a}, its spectrum is thrice that of~\(\fL_0=\tfrac12\int_{-1}^1u\,dy -u\) (definition~\eqref{eq:deffLn}).
Here it is straightforward to verify that the \(y\)-operator~\(\fL_0\) has: \begin{itemize}
\item one~\(0\) eigenvalue corresponding to eigenfunctions constant across the channel; and
\item  an `infinity' of eigenvalue~\(-1\) corresponding to all functions with zero average across the channel.
\end{itemize}

Then globally in space, with \(\de t\uv=\cL\uv+(\text{perturbation})\) at every \(X\in\XX\), and because of the `infinity' of the continuum~\XX, the linearised system has a `thrice-infinity' of the \(0\)~eigenvalue, and a `double-infinity' of eigenvalue~\(-1\).
Consequently, the theory of \cite{Aulbach2000} asserts: 
\begin{enumerate}
\item there exists a `\((3\infty)\)'-D slow manifold---the quasistationary~\eqref{eqqspdf}; 
\item which is  exponentially quickly attractive to all initial conditions, with transients roughly~\(e^{-t}\)---it is emergent; and
\item which we approximate by approximately solving the governing differential equations~\eqref{eq:losys}---done in encoded form by \cref{app}.
\end{enumerate}

We obtain a useful approximation to the global slow manifold by neglecting the `perturbing' remainder~\rv.
Because the remainder~\rv\ is the only coupling between different locales~\(X\) this approximation may be constructed independently at each and every cross-section~\(X\).
Further, because the Zappa system is homogeneous in space, the construction is identical at each and every \(X\in\XX\).
These two properties vastly simplify the construction of the attractive slow manifold.

Neglecting the coupling remainder~\(\rv\) gives the linear problem \(\de t\uv=\cL\uv\).
The approximate slow manifold is thus the zero eigenspace of~\cL.
We find the zero eigenspace via (generalised) eigenvectors.
With the (generalised) eigenvectors in the three columns of block-matrix~\cV, in essence we seek \(\uv(t)=\cV\Uv(t)\) such that \(d\Uv/dt=\cA\Uv\) for \(3\times3\) matrix~\cA\ having all the zero eigenvalues.
To be an eigenspace we need to solve \(\cL\cV=\cV\cA\)\,.
Now let's invoke previously established results \cite[\S2]{Roberts2013a}.
The linear operator~\(\cL\), defined in~\eqref{eq:losys}, has the same block Toeplitz structure as previously \cite[][(7) on p.1496]{Roberts2013a}.
Consequently \cite[Lemma~4]{Roberts2013a}, a basis for the zero eigenspace of~\cL\ is the collective columns of
\begin{align*}&
\cV=\begin{bmatrix} V_0&V_1&V_2
\\0&V_0&V_1\\0&0&V_0 \end{bmatrix},
&&\text{and further, }
\cA=\begin{bmatrix} 0&A_1&A_2
\\0&0&A_1\\0&0&0 \end{bmatrix}.
\end{align*}
The hierarchy of equations to solve for the components of these has been previously established \cite[Lemma~3]{Roberts2013a}: the hierarchy is essentially equivalent to the hierarchy one would solve if using the method of multiple scales, but the theoretical framework here is more powerful.
The upshot is that for Zappa dispersion, in which overlines denote cross-channel averages,
\begin{align}&
V_0=1,&& V_1=\overline v-v,&& V_2=2(\overline{v}^2-\overline{v^2}-\overline vv+v^2), 
\nonumber\\& 
&& A_1=-\overline v,&& A_2=\overline{(v-\overline v)^2}+\overline{v^2}\,.
\end{align}
In the specific case of \(v(y)=1-y^2\), these expressions reduce to the coefficients and polynomials of the slowly varying, slow manifold, model~\eqref{eqs:zappasm}.

So now we know that the evolution on the zero eigenspace, the approximate slow manifold, is \(\de t\Uv=\cA\Uv\): let's see how this translates into the macroscale \pde~\eqref{eq:zappact}.
Now, the first line of \(\de t\Uv=\cA\Uv\) is the \ode\ \(dU_0/dt=A_1 U_1+A_2U_2\).
Defining \(U_0=U(X,t):=\overline{u(X,y,t)}\), Proposition~6 of \cite{Roberts2013a} applies, and so generally~\(U(x,t)\) satisfies the macroscale effective advection-diffusion \pde~\eqref{eq:dcdtA}---a \pde\ that reduces to the specific~\eqref{eq:zappact} in the case \(v(y)=1-y^2\).

\subsection{Account for the coupling remainder}
Now we treat the exact `local' system \(d\uv/dt=\cL\uv+\rv\) as non-autonomously `forced' by coupling to all cross-sections in~\XX\ through the remainder \cite[aka Mori--Zwanzig transformation, e.g.,][]{Venturi2015}.
There are two justifications, both a simple and a rigorous, for being able to project such `forcing' onto the local model. 
First, simply, the rational projection of initial conditions for low-dimensional dynamical models leads to a cognate projection of any forcing \cite[\S7]{Roberts89b}. 
Second, alternatively and more rigorously, \cite{Aulbach2000} developed a general theory, that applies here, of centre manifolds for non-autonomous systems in suitable `infinite-D' state spaces: the theory establishes the existence and emergence of an `infinity-D' global centre manifold---a centre manifold whose construction \cite[Prop.~3.6]{Potzsche2006} happens to be symbolically identical at each \(X\in\XX\). 
Keep clear the contrasting points of view that contribute: on the one hand we consider the relatively low-dimensional system at each locale~\(X\) in space, a system that is weakly coupled to its neighbours; on the other-hand we consider the relatively high-dimensional system of all locales~\XX\ coupled together and then theory establishes global properties.

The upshot is that here we need to project the coupling remainder~\(\rv(t)\) onto each local slow manifold.  
Fortunately, the structure of the linear local dynamics~\eqref{eq:losys} is identical to that discussed by \cite{Roberts2013a}.
Hence, many of the results reported there apply here.
Linear algebra involving adjoint eigenvectors~\(Z_0\) and~\(\cW_n\) \cite[\(\fL^\dag_0Z_0=0\) and \(\cL^\dag\cW=\cW\cA\),][\S2.3]{Roberts2013a}, together with the history of the coupling remainder~\(e^{-t}\star\rv\), leads to the error formula~\eqref{eq:pderemain} \cite[equation~(23) from][]{Roberts2013a}.
Then the general macroscale advection-diffusion model~\eqref{eq:dcdtA} becomes exact with the error term~\(\rho\) included (here the error~\eqref{eq:pderemain} is third-order in spatial derivatives)
\begin{equation*}
\D tU= A_1\D xU+A_2\DD xU+\rho\,.
\end{equation*}
Then, simply, \emph{the macroscale effective advection-diffusion model \pde~\eqref{eq:dcdtA} is valid simply whenever and wherever the error term~\(\rho\) is acceptably small. 
There is: no~\(\epsilon\); no limit; no required scaling; no `balancing'; no ad hoc hierarchy of space-time variables.}

\section{Compact analysis, and connect to well-known methodology}
\label{seccacwm}

It is very tedious to perform all the algebraic machinations of \cref{secmkglm} on the Taylor series coefficients.  
Instead, we may compactify the analysis by defining the quadratic \emph{generating polynomial} \cite[\S3.1]{Roberts2013a}
\begin{equation}
\u(X,\x,y,t):=u_0(X,y,t)+\x u_1(X,Y,t)+\tfrac12\x^2 u_2(X,X,y,t)
\label{eqgpu}
\end{equation}
(or a higher-order polynomial if the analysis is to higher-order).
This generating polynomial then satisfies the exact differential equation~\eqref{eqgfut}. 
Consider \(\D t\u\), at \((X,\x,y,t)\), and substitute the equations~\eqref{eq:losys} for the Taylor coefficients at~\((X,y,t)\):
\begin{align}
\D t\u&=\D t{u_0}+\x\D t{u_1}+\tfrac12\x^2\D t{u_2}
\nonumber\\&
=\quad\fL_0u_0+\fL_1u_1+\fL_2u_2+r_0
\nonumber\\&\quad{}
+\fL_0\x u_1+\fL_1\x u_2 \qquad{}+\x r_1
\nonumber\\&\quad{}
+\fL_0\tfrac12\x^2u_2 \qquad\qquad{}+\tfrac12\x^2r_2
\nonumber\\&
=\fL_0\u+\fL_1\D\x\u+\fL_2\DD\x\u \ +\r
\nonumber\\
\implies \D t\u&
=\left[\fL_0{}+\fL_1\D\x{}+\fL_2\DD\x{}\right]\u+\r
\label{eqgfut}
\end{align}
for the generating polynomial of the coupling remainder, \(\r:=r_0+\x r_1+\tfrac12\x^2r_2\)\,.

Appropriate analysis of the \ide~\eqref{eqgfut} then reproduces the previous \eqref{secmkglm}.
But the algebra is done much more compactly as the separate components~\(u_0,u_1,u_2\) are all encompassed in the one generating polynomial~\u.
One important property of the analysis is that although we normally regard the derivative~\(\D\zeta{}\) as unbounded, in the analysis of \ide~\eqref{eqgfut} the space of functions is just that of quadratic polynomials in~\(\zeta\), and so here~\(\D\zeta{}\) is bounded, as well as possessing other nice properties.

Indeed, since we are only interested in the space of quadratic polynomials in~\(\zeta\), the analysis neglects any term~\Ord{\zeta^3}.  
Equivalently, we would work to `errors'~\Ord{\Dn\zeta 3{}}.
This view empowers us to organise the necessary algebra in a framework where we imagine~\(\D\zeta{}\) is `small'.
Note: in the methodology here \(\D\zeta{}\) is \emph{not assumed} small, as we track errors exactly in the remainder~\r, it is just that we may organise the algebra as if \(\D\zeta{}\) was small.
Such organisation then leads to the same hierarchy of problems as in \eqref{seclgsmt}, just more compactly.

\paragraph{Connect to extant methodology}
Since the notionally small~\(\D\zeta{}\) is effectively a small spatial derivative, we now connect to extant multiscale methods that a priori assume slow variations in space. 
That is, we now show that the non-remainder part of \ide~\eqref{eqgfut} appears in a conventional multiscale approximation of the governing microscale system~\eqref{eq:ukint}.

In conventional asymptotics we invoke restrictive scaling assumptions at the start.  
Here one would assume that the solution field~\(u(x,y,t)\) is slowly-varying in space~\(x\).
Then the argument goes that the field may be usefully written near any \(X\in\XX\) as the local Taylor quadratic approximation%
\footnote{I continue to conjecture that truncations to orders other than quadratic give corresponding analysis and results.  Ongoing research will elucidate.}
\begin{equation*}
u(\xi,y,t)\approx u|_{\xi=X}+\qt1u_\xi|_{\xi=X}+\qt2u_{\xi\xi}|_{\xi=X}\,.
\end{equation*}
Substituting into the nonlocal microscale~\eqref{eq:ukint} gives, at~\((X,y,t)\) and letting dashes\slash primes denote derivatives with respect to the first argument,
\begin{align}
\D tu&=\intYX k(X,\xi,y,\eta)u(\xi,\eta,t)\,d\xi\,d\eta
\nonumber\\&\approx \intYX k(X,\xi,y,\eta)\left[u|_{\xi=X}+\qt1u'|_{\xi=X}+\qt2u''|_{\xi=X}\right]d\xi\,d\eta
\nonumber\\&=\intYX  k(X,\xi,y,\eta)\,d\xi\,u(X,\eta,t)
+\int_{\XX}k(X,\xi,y,\eta)\qt1\,d\xi\,u'(X,\eta,t)
\nonumber\\&\qquad{}+\int_{\XX}k(X,\xi,y,\eta)\qt2\,d\xi\,u''(X,\eta,t)\ d\eta
\nonumber\\&=\int_{\YY}  k_0(X,y,\eta)u(X,\eta,t)
+k_1(X,y,\eta)u'(X,\eta,t)
\nonumber\\&\qquad{}
+k_2(X,y,\eta)u''(X,\eta,t)\,d\eta
\nonumber\\&=\fL_0u+\fL_1u'+\fL_2u''.
\label{equtfut}
\end{align}
Now the generating polynomial~\u, defined by~\eqref{eqgpu}, is such that \(u(X+\zeta,y,t)=\u(X,\zeta,y,t)+\Ord{\zeta^3}\).
Hence, rewriting the \emph{approximate} \pde~\eqref{equtfut} for~\(u(X+\zeta,y,t)\) at fixed~\(X\) gives precisely the \ide~\eqref{eqgfut} except that the remainder coupling~\r\ is omitted.
Consequently, extant multiscale methodologies continuing on from \pde
~\eqref{equtfut} generate equivalent results to that of~\cref{secmkglm}, but in a different framework---a framework without the error term~\eqref{eq:pderemain}.

Most extant multiscale analysis invokes, at the outset, balancing of scaling parameters, requires a small parameter, is only rigorous in the limit of infinite scale separation, and often invents heuristic multiple space-time variables.  
The approach developed herein connects with such analysis, but is considerably more flexible and, furthermore, justifies a more formal approach developed 30~years ago \cite[]{Roberts88a}, and implemented in \cref{app}. 
Further this approach derives the rigorous error expression~\eqref{eq:pderemain} at finite scale separation.

\section{Conclusion}

This article initiates a new multiscale modelling approach applied to a specific basic problem.
This article considers the scenario where the given physical problem~\eqref{eq:ukint} has non-local microscale interactions, such as inter-particle forces or dynamics on a lattice. 
Many extant mathematical methodologies derive, for such physical systems, an approximate macroscale \pde, such as the advection-diffusion~\eqref{eq:dcdtA}.
The novelty of our approach is that it derives a precise expression for the error of the macroscale approximate \pde, here~\eqref{eq:pderemain}.
Then, simply, and after microscale transients decay, the macroscale advection-diffusion \pde~\eqref{eq:dcdtA} is valid wherever and whenever the quantified error~\eqref{eq:pderemain} is acceptable.

Of course, in all such applications, we need the third moment of the microscale interaction kernel~\(k(x,\xi,y,\eta)\) to exist (see definition~\eqref{eq:du0dtr}) for the error analysis of \cref{secreld} to proceed and provide the error term.  
All moments exist for the Zappa problem, see~\eqref{eq:kern}.
If, in some application, the third moment does not exist, but the second moment does, then the advection-diffusion \pde~\eqref{eq:dcdtA} may be an appropriate macroscale model, but this work would not provide a quantifiable error.

Another important characteristic of our new approach is that the validity of the macroscale \pde\ is not confined by a limit `\(\epsilon\to0\)'---the approach holds for finite scale separation in the multiscale problem, in the large but finite domain~\XX.
Further, and in contrast to most extant methodologies, the approach here should generalise in further research to arbitrary order models just as it does when the microscale is expressed as \pde{}s \cite[]{Roberts2013a}.

The developed scenario here is that of linear nonlocal systems~\eqref{eq:ukint}.
However, key parts of the argument are justified with centre manifold theory \cite[e.g.]{Aulbach2000, Potzsche2006, Haragus2011, Roberts2014a}.
Consequently, further research should be able to show that cognate results hold for nonlinear microscale systems.
 
With further research, correct boundary conditions for the macroscale \pde{}s should be derivable by adapting earlier arguments to derive rigorous boundary conditions for approximate \pde{}s \cite[]{Roberts92c, Chen2016}.

Interesting applications of this novel approach would arise whenever there are microscale nonlocal interactions in the geometry of problems such as \cite[e.g.,][]{Roberts2014a}
dispersion in channels and pipes,
the lubrication flow of thin viscous fluids, 
shallow water approximations whether viscous or turbulent,
quasi-elastic beam theory,
long waves on heterogeneous lattices, and
\text{pattern evolution.}

\paragraph{Acknowledgements} This research was partly supported by the Australian Research Council with grant DP180100050.

\appendix
\section{Computer algebra derives macroscale PDE}
\label{app}
The following computer algebra derives the effective advection-diffusion \pde~\eqref{eq:zappact}, or any higher-order generalisation, for the microscale nonlocal Zappa system~\eqref{eq:muzap}. 
This code uses the free computer algebra package Reduce.%
\footnote{\url{http://www.reduce-algebra.com/}}
Analogous code will work for other computer algebra packages, and/or for cognate problems \cite[e.g.]{Roberts2014a}.

\begin{lstlisting}
% advection-diffusion PDE of Zappa transport in a channel
% AJR, 20 Jan 2017 -- 20 Jan 2020
on div; off allfac; on revpri; factor d,uu;

let d^5=>0; % truncate to this order of error
operator uu; depend uu,x,t; % uu(n):=df(uu,x,n)
let { df(uu(~n),x)=>uu(n+1), df(uu(~n),t)=>df(g,x,n)  };
operator mean; linear mean; % average across channel
let { mean(1,y)=>1, mean(y^~~p,y)=>(1+(-1)^p)/2/(p+1) };

% Preprocess nonlocal x-jumping: in essence finds the
% kernel integrals are (-v)^n
depend w,x; % dummy function for u(x)
% Taylor expand w(xi)=w(x+z) where z=xi-x
jmp:=for n:=0:deg((1+d)^99,d) sum d^n*df(w,x,n)*z^n/factorial(n)$
jmp:=int(exp(z/v)*jmp,z)$ % integrate exp((xi-x)/v)w(x)
% eval from z=-inf to 0 for the convolution
jmp:=sub(z=0,jmp/v)-w$ 

% iterate from quasi-equilibrium start
u:=uu(0)$ g:=0$
for it:=1:99 do begin
  res:=-df(u,t)+sub({w=u,v=1-y^2},jmp)+(-u+mean(u,y));
  write lengthres:=length(res);
  g:=g+(gd:=mean(res,y));
  u:=u+res-gd;
  if res=0 then write "Success: ",it:=it+10000;
end;
write "The resulting slow manifold and evolution is";
u:=u; duudt:=g;
end;
\end{lstlisting}

\end{document}